# The documentational approach to didactics

Luc Trouche, Ghislaine Gueudet, & Birgit Pepin[1]

**Abstract**

The documentational approach to didactics is a theory in mathematics education. Its first aim is to understand teachers' professional development by studying their interactions with the resources they use and design in/for their teaching. In this text we briefly describe the emergence of the approach, its theoretical sources, its main concepts and the associated methodology. We illustrate these aspects with examples from different research projects. This synthetic presentation is written for researchers, but also for non-specialists (e.g. master students) interested in a first discovery of the documentational approach.

## 1. Introduction

Mathematics teachers interact with curriculum and other resources in their daily work. Their work with resources includes selecting, modifying, and creating new resources, in-class and out-of-class. This creative work is termed *teacher documentation work*, and its outcome/s *teacher documentation.*

Typical curricular resources include text resources (e.g. textbooks, curricular guidelines, student worksheets); or digital curriculum resources (e.g. digital interactive textbooks). However, as there is now nearly unlimited access to resources on the web, teachers are often at a loss to choose the most didactically and qualitatively suitable resources for their mathematics teaching. Hence, the study of resources and mathematics teachers' interaction/work with those resources has become a prominent field of research (e.g. Pepin, Gueudet, & Trouche 2013), not least because curriculum reforms in many countries support the provision of reform-oriented curriculum materials that are seen to help teachers enact the new curriculum.

In theoretical terms the work of teachers with curriculum resources has been studied from many angles and theoretical perspectives (e.g. Remillard 2005; Pepin *et al.* 2013), for example in the Anglo/American research literature through the notion of 'enacted curriculum'. In continental Europe the notion of '*Didaktik*' is a common concept (e.g. Pepin *et al.* 2013). This entry describes, explains and illustrates the *Documentational Approach to Didactics* (DAD), which has its roots in French European Didactics.

## 2. Sources of the approach

The documentational approach to didactics (DAD) has been first introduced by Ghislaine Gueudet and Luc Trouche (Gueudet & Trouche 2009, Gueudet 2019), and has been developed further in joint work with Birgit Pepin (Gueudet, Pepin & Trouche 2012). DAD is originally steeped in the French didactics tradition in mathematics (Trouche 2016), where concepts such as *didactical situation*, *institutional constraint* and *scheme* are central. At the same time it also leans on socio-cultural theory, including notions such as *mediation* (Vygotsky 1978) as a constitutive of each cognitive process. Moreover, the approach has been developed due to the emerging digitalization of information and communication, which asks for new theoretical approaches.

The digitalization of information and communication and the development of Internet has indeed had strong consequences: ease of quick access to many resources and of communication with many people. This necessitated a complete metamorphosis of thinking and acting, particularly in education: new balances between *static* and *dynamic* resources, between *using* and *designing* teaching resources, between *individual* and *collective* work (Pepin, Choppin, Ruthven, & Sinclair 2017). Taking into account these phenomena, DAD proposed a change of paradigm by analyzing teachers' work through the lens of "resources" in and for teaching: the resources teachers prepare for (supporting) their classroom practices,

---

[1] This article is an updating of: Trouche, L., Gueudet, G., & Pepin, B. (2018). Documentational approach to didactics. In S. Lerman (Ed.), *Encyclopedia of Mathematics Education*. N.Y.: Springer. doi:10.1007/978-3-319-77487-9_100011-1

and that, in turn, get renewed by/in these practices.

In addition to the French didactics tradition, the authors drew their inspiration from several interrelated fields of research: the *field of technology use*, the field of *resources and curriculum design*, the field of *teacher professional learning/development*, the field of *information architecture,* and the field of *teacher professional development*.

In the field of technology use, the central underpinning for DAD was the *instrumental approach*. This theory has been developed by Rabardel (e.g. Vérillon & Rabardel 1995) in cognitive ergonomics, and then integrated into mathematics didactics (Guin, Ruthven, & Trouche 2005). It distinguishes between an *artefact*, available for a given user, and an *instrument*, which is developed by the user. Connected notions are those of *genesis*, of *instrumentation* and *instrumentalisation* – these are also essential component of DAD (see § 3). The development of the instrumental approach corresponded to a period where teachers were facing the integration of new *singular* tools (e.g. a calculator, a computer algebra software, a dynamic geometry system). It became clear that the instrumentation approach was not sufficient, as teachers were often surrounded (via Internet) by a profusion and a variety of *resources.*

This sensitivity to *resources* meets Adler's (2000) proposition of "think[ing] of a resource as the verb re-source, to source again or differently" (p. 207). Retaining this point of view, DAD took into consideration a wide spectrum of resources that have the potential to resource teacher activity (e.g. textbooks, digital resources, email exchanges with colleagues, student worksheets), resources *speaking to the teacher* (Remillard 2005) and supporting her/his engagement in teaching.

This broad view on resources lead to a broad view on teacher professional learning. As Ball, Hill and Bass (2005) stated in their study of 'mathematical knowledge for teaching', teaching is not reduced to the work in class, but also includes planning, evaluating, writing assessments, discussing with parents, amongst others. In DAD we consider the following: looking at teachers' work through their interactions with resources, and (following Cooney 1999) acknowledging that change of practice and change of professional knowledge and beliefs are connected (in a specific manner, as explained in Section 3).

Considering resources as the matter feeding teachers' work, a word was needed for naming what a teacher develops for a precise aim through his/her work with these resources. The word *document* was retained: it had already been used in the field of information architecture (Salaün 2012) for designing 'something bearing an intention', and dedicated to a given usage in a given context. This choice is the origin for the name of the approach, "documentational approach to didactics".

Finally, the easiness to communicate via the Internet lead to take into-account the emergence of a spectrum of various forms of teachers' collective work: networks, online association, communities more or less formal. Wenger's (1998) theory of *communities of practice,* and its concepts of *participation*, *negotiation*, and *reification* appeared as particularly fruitful for analyzing the design of teaching resources by collectives of teachers as a process of professional learning (Gueudet, Pepin, & Trouche 2013).

Once described, the sources of this theoretical approach, its structure and core concepts are presented in the following section.

## 3. The documentational approach to didactics - a holistic approach to teachers' work

In this section, the 'ingredients' of the DAD, and the processes involved, are described and explained. The following terms are defined: resources, documents, genesis, instrumentation, and instrumentalisation.

*Mathematics curriculum resources* are all the resources (e.g. digital interactive, non-digital/traditional text) that are developed and used by teachers and pupils in their interaction with mathematics in/for teaching and learning, inside and outside the classroom. This also includes digital curriculum resources (Pepin, Choppin, Ruthven, & Sinclair 2017), and Pepin and Gueudet (2018) make a distinction between digital curriculum resources and educational technology. They also distinguish between *material curriculum resources* (e.g. textbooks, digital curriculum resources, manipulatives and calculators), *social resources* (e.g. a conversation on the web/forum), and *cognitive resources* (e.g. frameworks/theoretical tools used to work with teachers). DAD has been mostly applied to teachers' work, but can also be used to study the work of teacher educators (e.g. Psycharis & Kalogeria, 2018), or students' interactions with resources (e.g. Kock & Pepin, 2018).

In terms of processes, during the interaction with a particular resource or sets of resources, teachers develop their particular *schemes of usage* with these resources (see section 4 below). These are likely to be different for different teachers, although they may use the same resource, depending on their dispositions and knowledge, for example. The outcome is the *document*, hence:

*Resources + scheme of usage = document*

The process of developing the document (including the teacher learning involved) has been coined *documentational genesis* (e.g., Gueudet & Trouche 2009).

Pepin, Gueudet, and Trouche (2013) have provided theoretical perspectives on 're-sourcing teachers' work and interactions', and the documentational approach is particularly pertinent to viewing the 'use' of resources as an interactive and potentially transformative process. This process works both ways: the affordances of the resource/s influence teachers' practice (the *instrumentation* process), as the teachers' dispositions and knowledge guide the choices and transformation processes between different resources (the *instrumentalisation* process) (Figure 1). Hence, the DAD emphasizes the dialectic nature of the teacher-resource interactions combining instrumentation and instrumentalisation (Vérillon & Rabardel, 1995). These processes include the design, re-design, or 'design-in-use' practices (where teachers change a document 'in the moment' and according to their instructional needs).

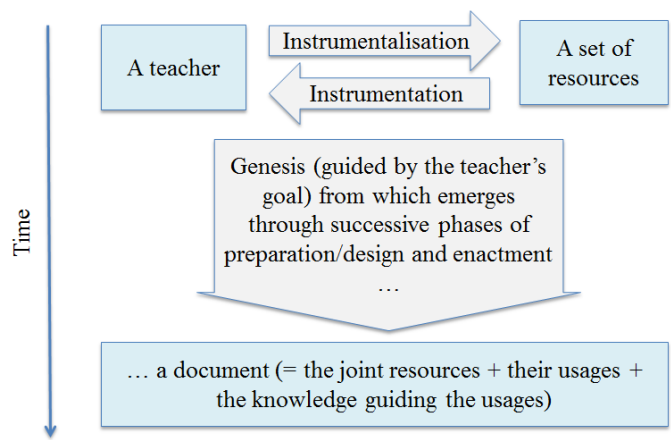

**Figure 1.** *A representation of a documentational genesis*

The DAD proposes a model of interactions between teachers and resources, which has implications for teachers' professional learning. Whilst there is an enormous amount of potentially suitable materials provided on the web, the web does not provide suitable support for relevant searching and selecting, which would be necessary if one wishes to search for particular (perhaps interactive) learning resources that combine with other (e.g. textbook) resources in their subtle epistemic or didactic features. In other words, what is provided for teachers is often "a pile of bricks", without giving guidance on how these bricks might be put together to develop a coherent student learning trajectory. Whether searching for tasks to supplement a given learning sequence, or planning learning paths through a flexible e-textbook, teachers will require professional support to help them develop *teacher design capacity* (Pepin, Gueudet, & Trouche 2017) - a mindfulness/sensitivity of mathematical and pedagogical aspects of learning resources, and the flexibility to use them (see Window 1). This is in line with the DAD, and Wang (2018) regards teacher design capacity as part of teachers' *documentational expertise.*

**Window 1. Resources at secondary school: the example of Vera's documentation work**

Vera is one of the many mathematics school teachers working with Sésamath in France (Gueudet *et al.* 2013; Pepin et al. 2017). Sésamath is an association of secondary school mathematics teachers in France, whose members have (since 2001) designed, and they freely offer, interactive e-textbooks on their website (http://www.sesamath.net/). Vera's documentation work is analyzed for a new lesson: it was the first time that Vera taught a grade 8 class on percentages, and she used a diversity of resources for this lesson, including Sésamath resources.

The analysis focuses on a lesson cycle: lesson preparation; enactment of the lesson; evaluation of student understanding; and reflection on her work. The choice of such a lesson cycle is in line with the

> ideas underpinning the DAD: the design was not restricted to the initial design of a given resource for teaching a particular content, but continued during the course of using the resource. The design work of Vera included, for example, the use of LaboMEP (a Sésamath tool) to propose different exercises to different students – this made her aware of the need to differentiate her teaching. LaboMEP also proposed variations of exercises with the same structure. Vera declared that it was a strong motivation for her to enhance her teaching, by mastering not only a set of familiar exercises but also those variations (with the same mathematical structure) related to a particular lesson.
>
> These results are interpreted as evidencing an enhancement of her didactical flexibility, in other words the development of her design capacity when selecting, utilizing and transforming existing curricular resources effectively, designing/creating new materials, for the purpose of effective mathematics instruction.

Teaching is often regarded as design. This is in line with Brown (2009) who explains that the interpretation of teaching as design, and the notion of teachers as designers, is fitting with a range of cognitive theories that "emphasize the vital partnership that exists between individuals and the tools they use to accomplish their goals. ... And it is not just the capacities of individuals that dictate human accomplishment, but also the affordances of the artifacts they use" (p. 19). Hence, Brown (2009) sees this relationship in the same way, as an interrelationship: the activity of "designing" is not only dependent on the teacher's competence, but it is an interrelationship between the teacher/s and the (curriculum) resources, the "teacher-tool relationship", that is at play here, and hence the affordances of the curriculum resources influence this relationship. This is in line with the DAD, emphasizing that any understanding of the teacher as designer must include a conscious/deliberate act of designing, of creating 'something new' (e.g., combining existing and novel elements) in order to reach a certain (didactical) aim (Pepin, Gueudet, & Trouche 2017) . This is provided and supported by *schemes of usage*, as defined in the following section.

### 4. Deepening the approach: schemes and systems

The concept of "scheme" (Vergnaud 1998) is central in the DAD. It is closely linked with the concept of "class of situations", which are in our context a set of professional situations corresponding to the same aim of the activity. For example, "managing the heterogeneity of the grade 8 class" is a class of situations for Vera (in the window given above). For a given class of situations, a subject (here a teacher) develops a stable organization of his/her activity, a scheme. A scheme has four components:

- The aim of the activity (characterizing the class of situations);
- Rules of action, of retrieving information and of control);
- Operational invariants, which are knowledge components of two (associated) kinds: *theorem-in-action* – a proposition considered as true- and *concept-in action* -a concept considered as relevant- (see example below);
- Possibilities of inferences, of adaptation to the variety of situations.

Over the course of his/her activity, the teacher can enrich his/her schemes, integrating new rules of actions, or s/he can develop new schemes: the scheme offers in fact a model for analyzing learning. In DAD, the schemes that are considered are schemes of usage for a given resource (or set of resources). The resources and the scheme make up a document (as summarized by the equation provided in the previous section).

The set formed by all the resources used by the teacher is named his/her *resource system*. These resources are associated with schemes of usage, forming documents (the same resource can intervene in several documents). The documents developed by a teacher also form a system, called the *document system* of the teacher. Its structure follows the structure of the class of situations composing the professional activity of the teacher (according to the different aims of his/her activity).

When teachers share their documentation work, for example in a group preparing lessons collectively, they may also develop a shared resource system (Trouche, Gueudet & Pepin 2019). Nevertheless, the different members of the group can develop different schemes for the same resource, resulting in different documents (Pepin & Gueudet 2020).

Window 2 presents a case of a resource at primary school, in order to illustrate/exemplify operational invariants, resource systems and document systems.

**Window 2. Resources at primary school: the example of the virtual abacus**

The virtual abacus (Figure 2) is a free software[2] developed in France by Sésamath, an association of mathematics teachers designing online resources (see Window 1).

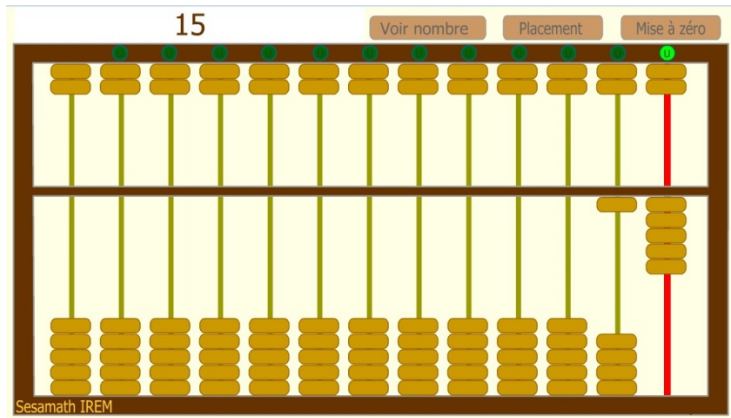

**Figure 2.** *The virtual abacus*

The Chinese abacus is separated in two parts by a central bar, called "the reading bar": only the beads on this bar are considered as "activated". There are two kinds of beads: 5-unit beads (two of them) and 1-unit beads (five of them). The Chinese abacus comprises 13 vertical rods. Each rod corresponds to a rank of the place-value system: units, tens, hundreds, etc. (from the right to the left). There are several possibilities to display the same number on the Chinese abacus: for example 15 is represented on the abacus above using 6 beads (one-unit beads, one on the tenth rod and five on the unit rod); it can also be represented using only two beads, by replacing the five one-unit beads on the right by a five-units bead on the same rod.

Carlos is an experimented primary school teacher who has been "shadowed" over three years (Poisard, Bueno-Ravel, & Gueudet 2011). He has decided to use the abacus for his teaching of numbers with his grade 3 class. He has seen in his mathematics textbook an activity with the abacus, but did not want to use it in class before discovering the virtual abacus. He started with material abacii, the students manipulated them and formulated hypotheses on the way it works. Then they worked on the virtual abacus, and wrote instructions for its use. After that Carlos proposed exercises: inscribing a given number on the abacus, reading a number inscribed on the abacus. For the final assessment such exercises were given on paper, to avoid a strategy of trial and error possible on the software.

Carlos developed several documents incorporating the virtual abacus and other associated resources (Poisard, Bueno-Ravel, & Gueudet ibid). For the aim "Discover how the abacus works", he used both material abaci and the virtual abacus, and asked the students to compose posters. It was important for him to let his students discover by themselves the principles of the abacus. This corresponds to an operational invariant: a theorem-in-action like "the students must discover by themselves a far as possible the new tools they meet". It is also linked to the associated concept-in-action: "self-discovery", and both were developed by Carlos before meeting the virtual abacus. Another operational invariant intervened in his choices: "it is important for students in grade 3 to manipulate material resources". The scheme for the aim "discover how the abacus works" comprises these operational invariants and associated rules of actions: "propose material abaci to the students to afford manipulation"; "propose the virtual abacus to the students to allow them to check which number is displayed".

Along his work with the abacus, Carlos observed that for the exercise: "display a given number on the abacus", if the students use the virtual abacus, they develop trial and error strategies using the button "display number". Hence he decided to do a final assessment on paper. He developed a new document, for the aim: "teach the students how to display a number on the abacus". This document incorporates the virtual abacus, but also abaci drawn on paper, and an operational invariant like "on the virtual abacus students can use trial and error strategies".

---

[2] The software and other resources can be found at http://seminaire-education.espe-bretagne.fr/?page_id=611

Carlos is an experienced primary school teacher. For the teaching of numbers in grade 3, he had for many years developed a (sub-)system of resources and documents. The abacus has been incorporated in these resources, and in new documents. Some of these documents correspond to aims directly linked with the abacus, like "discover how the abacus works", "teach the students how to display a number on the abacus". In other documents, the abacus did not appear in the aim, but was nevertheless used for this aim. For example, for the aim "teach the principles of the base ten place-value system", he used the abacus to evidence the "grouping and exchanging" principles (like grouping two five-beads on a rod and exchanging it with a one-unit bead on the next rod). Other resources in his resource system intervened in these documents, like posters written by the students. Some of these resources have been decisive in his choice to use the abacus: the textbook in particular, which offered him a first opportunity to meet a possible use of the abacus in class.

The whole documents system of a teacher comprises many sub-systems with their own structure, linked for example with a particular mathematical content or a specific kind of activity: sub-system for geometry, sub-system for assessment. It can be described at different levels, ranging from a very general view of the activity to a very specific focus on a given mathematical content. For research in mathematics education, the more specific levels, taking into account the mathematics content (the aim could be "assess the students' skills on percentages in grade 8"), are more informative about teachers' interactions with resources and their consequences. DAD claims that analyzing teachers' documentation work requires specific methodology, which is the purpose of the following section.

## 5. Reflective investigation: a developing methodological construct

This section presents the research design typically linked to the DAD; and the principles grounding this design; then describes one tool illustrating these principles; finally presents some issues that the reflective investigation methodology is facing. Analyzing teachers' activity through their documentation work requires to take-into-account the following: the variety of resources feeding, and produced by, this work; the variety of interactions (collective, institutional as well as social) influencing this work; the time for developing documentational geneses. These epistemological considerations lead DAD to develop a specific methodology, named *reflective investigation of teachers' documentation work*.

This methodology gives a major role to teachers themselves, and it is underpinned by five main principles[3]:

- The principle of *broad collection of the material resources* used and produced in the course of the documentation work;
- The principle of *long-term follow-up*. Geneses are ongoing processes and schemes develop over long periods of time;
- The principle of *in- and out-of-class follow-up*. The classroom is an important place where the teaching elaborated is implemented, bringing adaptations, revisions and improvisations. However, an important part of teachers' work takes place outside the classroom: at school (e.g. in staff rooms), at home, in teacher development centers/programs;
- The principle of *reflective follow-up* of the documentation work;
- The principle of *confronting the teacher' views on her documentation work, and the materiality of this work* (materiality coming for example from the collection of material resources; from the teacher's practices in her classrooms).

The active involvement of the teacher is a practical necessity, as s/he is the one having access to his/her documentation work (beyond the direct observation of the researcher). It also yields a reflective stance leading the teacher to an introspective attitude, sometimes making visible what could be hidden resources, or hidden links within his/her resource system. The principles and the purpose of this data collection spread over a long period of time must be presented clearly by the researcher to the teacher from the beginning. This necessity leads Sabra (2016) to propose the notion of *methodological contract* linking a teacher and a researcher following his/her documentation work.

---

[3] In (Gueudet et Trouche 2012, p. 27), we proposed only the four first principles. We add here a fifth one, which progressively emerged from the studies developed in the frame of DAD.s

Based on these principles, selected data collection strategies and tools were designed, adapted to the various contexts and research questions. For example, a tool seemingly fruitful is the *schematic representation of a teacher's resource system* (SRRS, see Figure 3). The teacher is asked to draw a map of her resources, evidencing the resources s/he had identified and appropriated, from which repositories, and for which purpose (e.g. Pepin, Xu, Trouche, Wang, 2016).

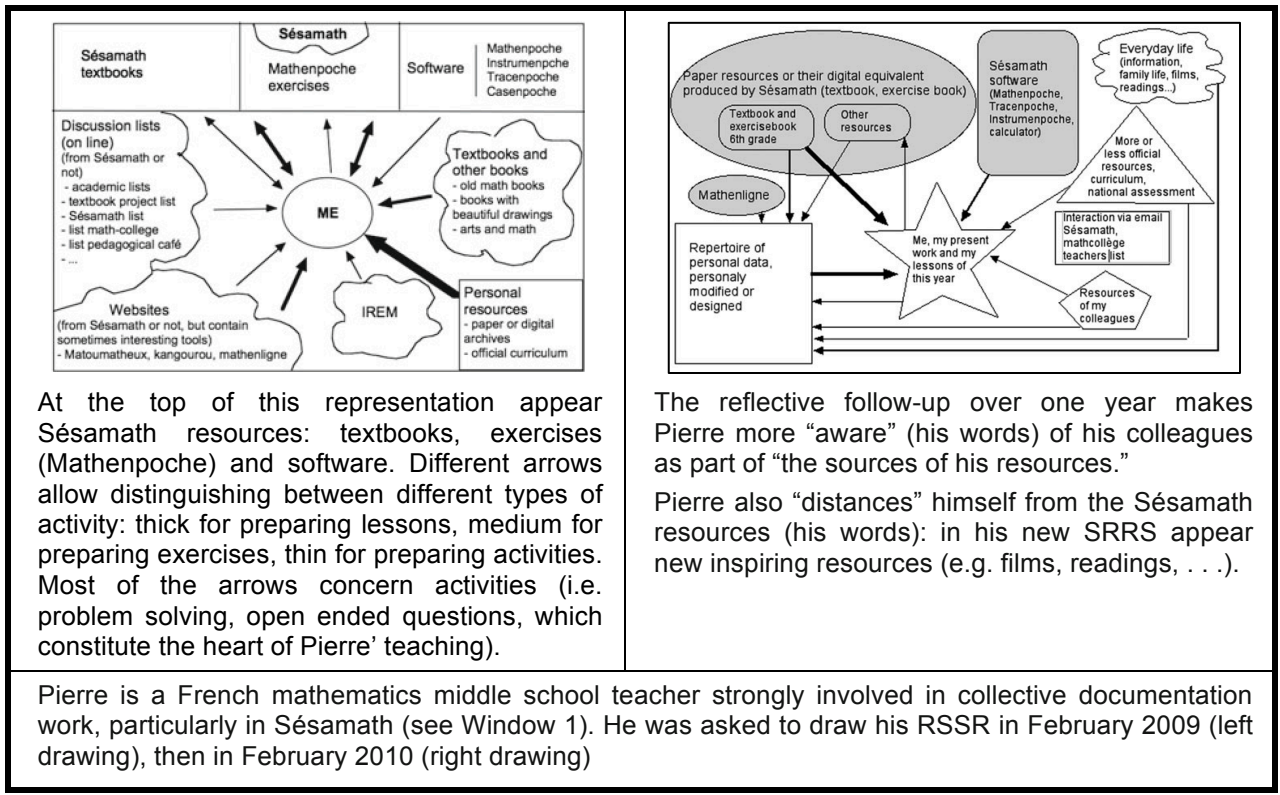

At the top of this representation appear Sésamath resources: textbooks, exercises (Mathenpoche) and software. Different arrows allow distinguishing between different types of activity: thick for preparing lessons, medium for preparing exercises, thin for preparing activities. Most of the arrows concern activities (i.e. problem solving, open ended questions, which constitute the heart of Pierre' teaching).

The reflective follow-up over one year makes Pierre more "aware" (his words) of his colleagues as part of "the sources of his resources."

Pierre also "distances" himself from the Sésamath resources (his words): in his new SRRS appear new inspiring resources (e.g. films, readings, . . .).

Pierre is a French mathematics middle school teacher strongly involved in collective documentation work, particularly in Sésamath (see Window 1). He was asked to draw his RSSR in February 2009 (left drawing), then in February 2010 (right drawing)

**Figure 3.** *A teacher's SRRS (Gueudet, Pepin & Trouche 2012, p. 314 & 318)*

Since the beginning of DAD, this tool has been developed in several directions:

- Hammoud (2012), working in chemistry education, proposed an approach for analyzing SRRS as mind maps; she also used the SRRS for asking teachers to describe their interactions with colleagues, or within different collectives;
- Rocha (2018) renamed SRRS as 'Reflective mapping of teacher resource system' (RMRS), for two reasons: emphasizing the role of reflectivity; and denoting (with 'mapping') a process of progressive exploration of an unknown territory (for the researcher, but to some extent also for the teacher him/herself). In her research she asked a teacher to draw such RMRSs at different moments, to investigate to what extent particular 'teaching' and reform developments (linked to particular resources) result in a « new » RMRS. In addition, following the fifth principle (see above), she confronted the teacher with a given RMRS, made by the teacher, and developed what she named 'Inferred mapping of teacher resource system' (IMRS) made by the researcher herself.

Beyond this specific tool, new methodological developments of the reflective investigation have occurred, to access, as far as possible, the 'true' teacher documentation work at the time, when it happens. It means observing a teacher during episodes of interaction with resources, not only in public episodes (e.g. in classrooms) but also during more 'intimate' episodes: for example, preparing a progression for the year; preparing a lesson; revising it. A video-recorded follow-up of these successive episodes has been conducted (Bellemain & Trouche 2016), giving access to gestures and words of a teacher in the course of her documentation work (allowing, for example, to infer elements of schemes). This kind of follow-up raises some difficulties, amongst them the following:

- Under which 'natural' conditions could a teacher work alone with resources, and at the same time describe the rationale of her activity? Wang (2018) introduces the notion of teacher's *documentation-workmate*, meaning a teacher regularly sharing the documentation work with a colleague (the documentation-workmate). The follow-up of the pair of teachers working together gave access to their mutual explanations, and to aspects of the teacher's knowledge guiding her documentation work (Trouche et al. 2019);
- How is it possible to store (for analysis purposes) the heterogeneous and numerous data resulting from the follow-up of teachers' documentation work? This issue is addressed in the frame of the AnA.doc project (Alturkmani *et al.* 2019), from the development of a platform prototype. It allows storing videos as well as resources intervening or produced by teachers' documentation work, for subsequent analysis and sharing of both the data and their analyses within a research community.

Analyzing teachers' documentation work is a complex process. The five principles of the reflective investigation methodology provide guidelines for methodological choices. Researchers using these principles need to make motivated choices for limiting the abundance of data: choices of critical moments for teachers' documentation work (see the notion of *documentational incidents*, Sabra 2016); or critical resources for a teacher's resource system (see the notion of *pivotal resource*, Gueudet 2017). This work is in progress, and the methodological and conceptual issues need to be stronger interrelated.

## 6. Perspectives for further evolutions

The documentational approach to didactics is a recent theoretical framework in mathematics education. While the concepts presented in section 2 and 3 (resources, documents, documentational geneses) are now well established, more elaborated concepts like 'resource systems' and 'document systems' are still evolving. For example, in the case of secondary school mathematics teachers, is it possible to observe different types of structures of their document systems that can be associated with particular professional profiles of these teachers? Do some groups of teachers (communities of practice, in particular) share collective document systems (e.g. Pepin & Gueudet 2020)?

Moreover, the fields of application of the DAD have evolved over the years. In terms of education levels they now range from preschool (Besnier & Gueudet 2016) to university (Gueudet 2017; Kock & Pepin, 2018), and also include the work of teacher educators with resources (Psycharis & Kalogeria 2018). In terms of disciplines, the documentational approach has been used in experimental sciences like physics and chemistry (Hammoud 2012), and also in language education (Quéré 2019). What are the specificities of teachers' documentation work and documentation systems in these new contexts?

The use of DAD in various social and cultural contexts (e.g. Brazil, China, Lebanon, Norway, Senegal) also leads to questions about the 'cultural coloring' of teachers' work with resources: for example, different naming systems used by teachers in their daily documentation work; or the perceptions and practices related to 'teacher design' vary amongst cultural spheres (e.g. Pepin, Artigue, Gitirana, Miyakawa, Ruthven, & Xu 2019). This diversity is the result of historical, social and cultural contexts in which teachers' work takes place. Research investigating these questions could lead to a better understanding of a more nuanced view of the nature of teachers' interactions with resources and to a deepening of the DAD related concepts.

In a recent study the work of students with resources has been investigated (e.g. Kock & Pepin 2018). Other studies (e.g. Gueudet & Pepin 2018; Trouche, Gitirana, Miyakawa, Pepin, & Wang 2019) have suggested possible links with other theories. The theory of didactical situations (Brousseau 1998), for example, introduces the concept of *milieu* which includes all the objects with which the individual student interacts in a mathematical situation. These objects can be considered as resources. What are the consequences of such theoretical links?

As the range of teaching and learning phenomena studied with this approach increases, the concepts and methods also develop. This is evident in numerous communications, publications and specialized conferences (e.g. Re(s)sources 2018 International Conference, see Gitirana *et al.* 2018 and Trouche, Gueudet & Pepin 2019) using DAD. Such events and publications also evidence that there are 'missing resources' to be developed (Trouche 2019). These studies and activities constitute milestones in the ongoing development of a living theoretical framework.